\documentclass[11pt]{article}
\usepackage[english]{babel}
\usepackage[a4paper]{geometry}
\usepackage{amsmath,amsfonts,amssymb,amsthm}
\usepackage{graphicx}
\usepackage[colorlinks=true, allcolors=blue]{hyperref}
\usepackage{array}
\usepackage{breqn}
\usepackage{comment}
\usepackage{enumitem}
\usepackage{color}
\usepackage{tabularx}
\usepackage{csquotes}
\usepackage{authblk}
\usepackage{caption}
\usepackage{multicol}

\newcolumntype{Y}{>{\centering\arraybackslash}X}

\usepackage[maxbibnames=99,backend=bibtex,style=numeric]{biblatex}
\usepackage{graphicx}
\graphicspath{ {./images/} }

\newtheorem{definition}{Definition}[section]
\newtheorem{theorem}[definition]{Theorem}
\newtheorem{lemma}[definition]{Lemma}

\newcommand{\ceilhalf}{f_0}

\title{Sparse metric hypergraphs}
\author[1]{Vašek Chvátal\thanks{Supported by H2020-MSCA-RISE project CoSP- GA No. 823748, and by NSERC (Natural Sciences and Engineering Research Council of Canada) grant RGPIN/5599-2014. E-mail: vchvatal9@gmail.com}}

\author[2]{Guillermo A. Gamboa Q.\thanks{Supported by project 22-17398S (Flows and cycles in graphs on surfaces) of Czech Science Foundation. E-mail: gamboa@iuuk.mff.cuni.cz}}

\author[3]{Ida Kantor\thanks{Supported by project 22-19073S of Czech Science Foundation and by Charles University project ˇ
UNCE/SCI/004. E-mail: ida@kam.mff.cuni.cz}}

\affil[1,2,3]{Charles University, Prague, Czechia}
\affil[1]{Concordia University, Montreal, Canada}
\date{}

\begin{document}
\maketitle

\begin{abstract}
Given a metric space $(X, \rho)$, we say $y$ is \textit{between} $x$ and $z$ if $\rho(x,z) = \rho(x,y) + \rho(y,z)$. A metric space gives rise to a 3-uniform hypergraph that has as hyperedges those triples $\{ x,y,z \}$ where $y$ is between $x$ and $z$. Such hypergraphs are called \textit{metric} and understanding them is key to the study of metric spaces. In this paper, we prove that hypergraphs where small subsets of vertices induce few edges are metric. Additionally, we adapt the notion of sparsity with respect to monotone increasing functions, classify hypergraphs that exhibit this version of sparsity and prove that they are metric.
\end{abstract}

\section{Metric Hypergraphs}

In 2003, Chen and Chv\'{a}tal initiated an effort to prove the analogues of various theorems from Euclidean geometry within the realm of finite metric spaces. Chen managed to prove the analogue of Sylvester-Gallai theorem in this setting \cite{chenfano}, and there are many interesting partial results by several groups of researchers leading towards the analogue of de Bruijn--Erd\H{o}s theorem \cite{originalconj,mainpaper,numberoflines,ik}

A central notion in these theorems is that of {\em betweenness}: if $x,y,z$ are three distinct points of a metric space $M=(X,\rho)$, we say that $y$ is {\em between} $x$ and $z$ if 
\[\rho(x,z)=\rho(x,y)+\rho(y,z).\]
 If one of $x,y,z$ is between the other two, we say that the triple $\{x,y,z\}$ is {\em collinear}. 

While studying the questions mentioned above, it turned out that a particularly useful tool is the 3-uniform {\em hypergraph $H_M$ associated with the  metric space $M$}: 
\begin{itemize}
    \item the vertices of $H_M$ are the points of $M$, 
    \item the hyperedges are all triples $\{x,y,z\}$ that are collinear in $M$. 
\end{itemize}
Each metric space has a hypergraph associated with it, but not every 3-uniform hypergraph is associated with a metric space. We say that a 3-uniform hypergraph $H$ is {\em metric} if there is a metric space $M$ such that $H_M= H$. There are several examples of hypergraphs that are not metric. For each $n$, a family of $\Omega(2^n)$ non-metric hypergraphs on $n$ vertices was constructed by Chv\'{a}tal and Kantor in \cite{chknonmetri}. It follows from the result of Chen \cite{chenfano} that Steiner triple systems of order greater than 3 are not metric. 

It would be very useful to have an efficient way of recognizing metric hypergraphs. One characteristic that seems to have a lot to do with whether a hypergraph is metric is some notion of ``density"---sparser hypergraphs seem more likely to be metric than denser ones. In particular, Chvátal and Kantor \cite{chknote} proved that any noncomplete 3-uniform hypergraph with $n$ vertices and at least $\binom{n}{3} - n+5$ hyperedges is nonmetric. As they also described some metric hypergraphs with $\binom{n}{3}-n+4$ hyperedges, this bound is tight.

Inspired by this, one might ask if there is some function $f(n)$ such that all hypergraphs on $n$ vertices with at most $f(n)$ hyperedges are metric. However, an answer to this question is trivial---if we allow at least 7 edges, then we may construct a nonmetric hypergraph by adding isolated vertices to Fano plane (which, being a Steiner triple system, is nonmetric). This trivial example contains a subhypergraph on a very few vertices, with many edges. It seems reasonable to suspect that forbidding such small dense subgraphs might guarantee metricity.

Our first result is in this direction. 

\begin{definition}
\label{defsixtwo}
    We say that a 3-uniform hypergraph $H$ is $(k,\ell)$-sparse if for any k-element subset $X$ of vertices, there are at most $\ell$ hyperedges with all vertices contained in $X$.
\end{definition}

\begin{theorem}
\label{thm:sixtwo}
If $H$ is a $(6,2)$-sparse 3-uniform hypergraph, then $H$ is metric.  
\end{theorem}

Next natural candidates for metricity are $(5,2)$-sparse and $(6,3)$-sparse hypergraphs. Unfortunately, these conditions do not guarantee metricity, since 
the Fano plane is $(5,2)$-sparse, and there are several constructions of Steiner triple systems that are $(6,3)$-sparse (these are sometimes called {\em anti-Pasch} Steiner triple systems), see \cite{MR1762019}. 

Steiner triple systems may be locally sparse in some sense---in fact, Erd\H{o}s conjectured in 1976 (\cite{MR0465878}) that for every $k$, if $n$ is a large enough integer that satisfies the natural divisibility condition for existence of Steiner triple systems, then there is a $(k,k-3)$--sparse Steiner triple system of order $n$. This was proved by and Glock, K\"{u}hn, Lo, and Osthus in 2020 (\cite{MR4121151}). But Steiner triple systems have larger subsets with many hyperedges, namely subsets of $k$ vertices with $\Omega(k^2)$ hyperedges (for $k$ close to the number of vertices). So it might be that $(6,3)$-sparse hypergraphs are metric if we impose additional conditions on larger subsets of vertices. This motivates the following definition, and our second result.

\begin{definition}
\label{def:fsparse}
Let $f(x)$ be a monotonically increasing function. We say that $H$ is $f$-sparse if for every $k \in \{4, 5, \ldots, n\}$, any $k$-set $X \subset V$ induces at most $f(k)$ hyperedges. Our second result concerns the function $f_0$ defined as $f_0(k)=\lceil \frac{k}{2}\rceil$.
\end{definition}

\begin{theorem}\label{thm:fsparse}
	If $H$ is a $f_0$-sparse 3-uniform hypergraph, then $H$ is metric. 
\end{theorem}

{\noindent\bf Problem.} {\em What about other functions? For instance, are all $f$-sparse hypergraphs metric, if $f(k)=2k/3$, or $f(k)=k-4$?}
\vspace{3mm}

{\noindent\bf Problem.} {\em Are there (infinite families of) nonmetric hypergraphs with ``few edges" that do not contain Steiner triple systems as subhypergraphs?}

\section{(6,2)-sparse hypergraphs are metric}

\begin{theorem}\label{sixtwoNoGoat}
    If $H$ is a $(6,2)$-sparse $3$-uniform hypergraph where two hyperedges share at most one vertex, then $H$ is metric.
 \end{theorem}

 \begin{proof}
     For every hyperedge $e \in E$, we pick an arbitrary vertex and label it $x_e$. Note that a vertex $x$ can receive 
    several labels. Let $\rho_0: V^2 \rightarrow \{ 0,1,1.5,2 \}$ be the map defined as
    \begin{equation*}
        \rho_0 (x,y) =
        \begin{cases}
            0 &\text{if } x=y, \\
            1 &\text{if there is a hyperedge $e$ such that $\{ x,y \} \subset e$ and $x_e \in \{x,y\}$},  \\
            2 &\text{if there is a hyperedge $e$ such that $\{ x,y \} \subset e$ and $x_e \notin \{x,y\}$}, \\
            1.5 &\text{otherwise.}
        \end{cases}
    \end{equation*}
    Since two hyperedges share at most one vertex, $\rho_0$ is well-defined. It is also clear that $\rho$ is a metric. Consider the metric space $M=(V,\rho_0)$. 
    
    If $\{ x_e, y, z \}$ is a hyperedge of $H$, then $ \{ x_e, y ,z \}$ is a collinear triple in $M$ since $\rho_0(y,z) = \rho_0(y, x_e) + \rho_0(x_e, z)$. 
    
    On the other hand, suppose that $\{x,y,z\}$ is collinear in $M$, with $x$ between $y$ and $z$. Then $\rho_0(x,y)=\rho_0(y,z)=1$ and $\rho_0(x,z)=2$. Since $\rho_0(x,y)=1$, the vertices $x,y$ are together in some hyperedge $e_{xy}$ of $H$. Similarly, $y,z$ are together in some hyperedge $e_{yz}$, and $x,z$ are together in $e_{xz}$. If these are three different hyperedges, then $H$ is not $(6,2)$-sparse. The only other option is $e_{xy}=e_{yz}=e_{xz}=\{x,y,z\}$.
 \end{proof}

In a complete analogy with graphs, a hypergraph $H=(V,E)$ is {\em connected} if we cannot partition the vertex sets into disjoint subsets $V_1$ and $V_2$ such that every hyperedge is a subset of $V_1$ or of $V_2$. The following lemma shows that we can restrict our attention to connected hypergraphs.

\begin{lemma}
    \label{connected}
    If $H_1=(V_1,E_1)$ and $H_2=(V_2,E_2)$ are metric hypergraphs such that $V_1\cap V_2=\emptyset$, then $H=(V_1\cup V_2, E_1\cup E_2)$ is a metric hypergraph as well.
\end{lemma}
    
\begin{proof}
    For $i\in \{1,2\}$, let $M_i=(V_i,\rho_i)$ be a metric space such that $H_{M_i}=H_i$. Let $V=V_1\cup V_2$ and let $m$ be a number such that for all $x,y\in V_1$, $\rho_1(x,y)<m$, and for all $x,y\in V_2$, $\rho_2(x,y)<m$. Define $\rho:(V\times V)\to \mathbb{R}$ as follows: 
        \begin{equation*}
        \rho (x,y) =
        \begin{cases}
            \rho_1(x,y) &\text{if } x,y\in V_1, \\
            \rho_2(x,y) &\text{if } x,y\in V_2, \\
            m &\text{if  $x\in V_1$ and $y\in V_2$ or vice versa} 
        \end{cases}
    \end{equation*}
    Then $M=(V,\rho)$ is a metric space such that $H_M=H$. 
\end{proof}

\begin{lemma}
\label{lemma1}
    If $H=(V,E)$ is a connected $(6,2)$-sparse 3-uniform hypergraph with $|V|\geq 6$, then two distinct hyperedges contain at most one common vertex. 
\end{lemma}

We are now ready to prove Theorem~\ref{thm:sixtwo}.

\begin{proof}[Proof of Theorem~\ref{thm:sixtwo}]
If $H$ has at most $5$ vertices, then $H$ has either only one edge, or two edges sharing one or two vertices. In each case it is easy to find a corresponding metric space.

If $H$ has at least $6$ vertices, we may assume that $H$ is connected, otherwise we find a metric space for each component separately and use Lemma~\ref{connected}. By Lemma~\ref{lemma1}, two hyperedges share at most one vertex, so the conclusion follows by Theorem~\ref{sixtwoNoGoat}.   
\end{proof}

\section{Sparsity in terms of monotonically increasing functions}

In this section, we will focus on $\ceilhalf$-sparse hypergraphs with $f_0$ given by $f_0(k) = \lceil \frac{k}{2} \rceil$.

Let $H=(V,E)$ be a hypergraph and $X \subseteq V$. The \textit{restriction} of $H$ to $X$, denoted by $H[X]$, is the hypergraph with vertex set $X$, and hyperedge set $E[X] := \{ e \in E \vert e \subseteq X\}.$ The hypergraph $(V, E \setminus E[X])$ will be denoted by $H \setminus E[X]$. We say that the hyperedges in $E[X]$ are \textit{induced} by $X$.

In order to prove our main result, we characterize $\ceilhalf$-sparse hypergraphs. The proof is somewhat technical, so we omit it here in an extended abstract.

\begin{figure}[ht]
\begin{multicols}{6}
\includegraphics[width=\linewidth]{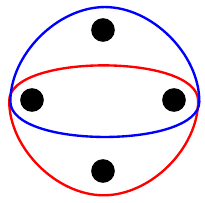}\par
\includegraphics[width=\linewidth]{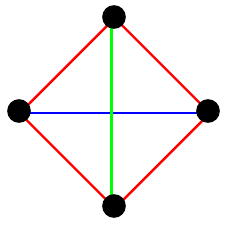}\par
\includegraphics[width=\linewidth]{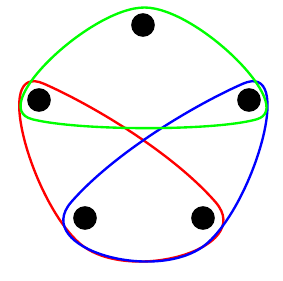}\par
\includegraphics[width=\linewidth]{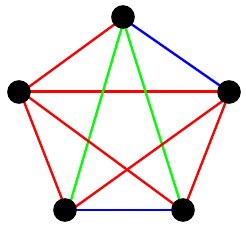}\par
\includegraphics[width=\linewidth]{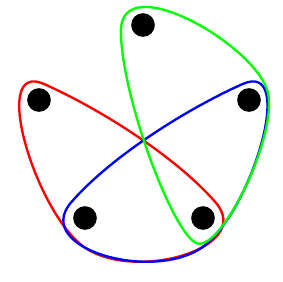}\par
\includegraphics[width=\linewidth]{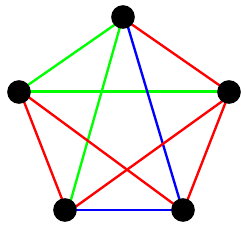}\par
\end{multicols}
\begin{multicols}{3}
    \begin{center}
        $H_4$ and $M_4$
    \end{center}\par
    \begin{center}
        $H_5^1$ and $M_5^1$
    \end{center}
    \begin{center}
        $H_5^2$ and $M_5^2$
    \end{center}
\end{multicols}
\begin{multicols}{6}
\includegraphics[width=\linewidth]{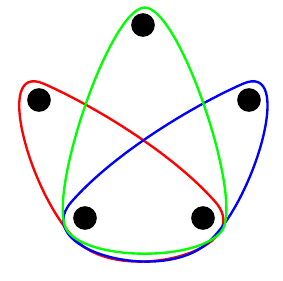}\par
\includegraphics[width=\linewidth]{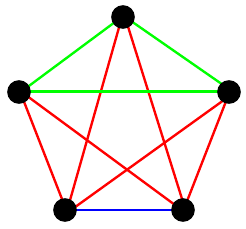}\par
\includegraphics[width=\linewidth]{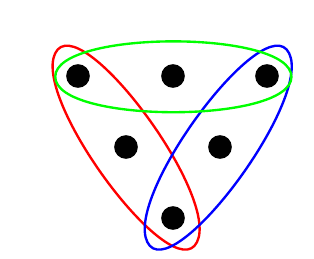}\par
\includegraphics[width=\linewidth]{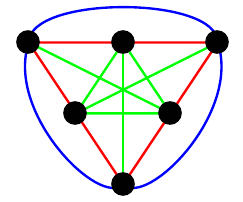}\par
\includegraphics[width=\linewidth]{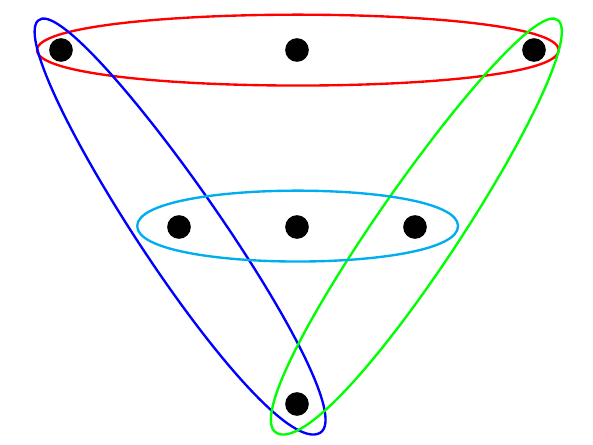}\par
\includegraphics[width=\linewidth]{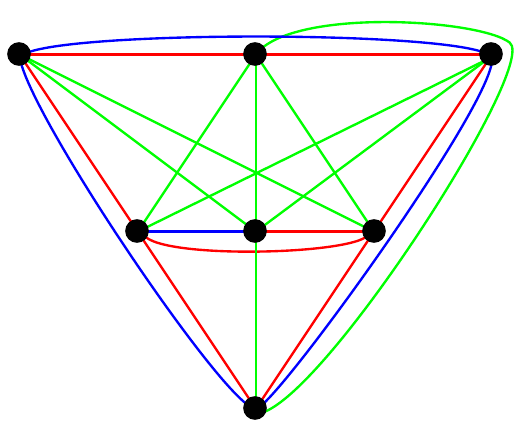}\par
\end{multicols}
\begin{multicols}{3}
    \begin{center}
        $H_5^3$ and $M_5^3$
    \end{center}\par
    \begin{center}
        $H_6$ and $M_6$
    \end{center}
    \begin{center}
        $H_7^1$ and $M_7^1$
    \end{center}
\end{multicols}
\begin{multicols}{4}
\includegraphics[width=2cm]{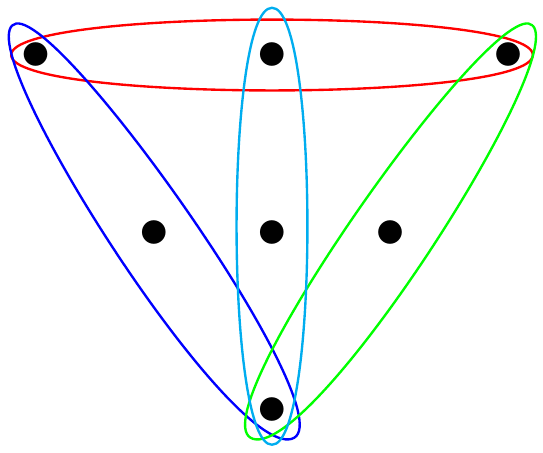}\par
\includegraphics[width=2cm]{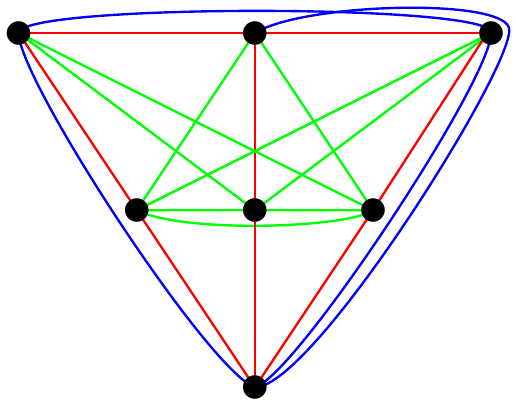}\par
\includegraphics[width=2cm]{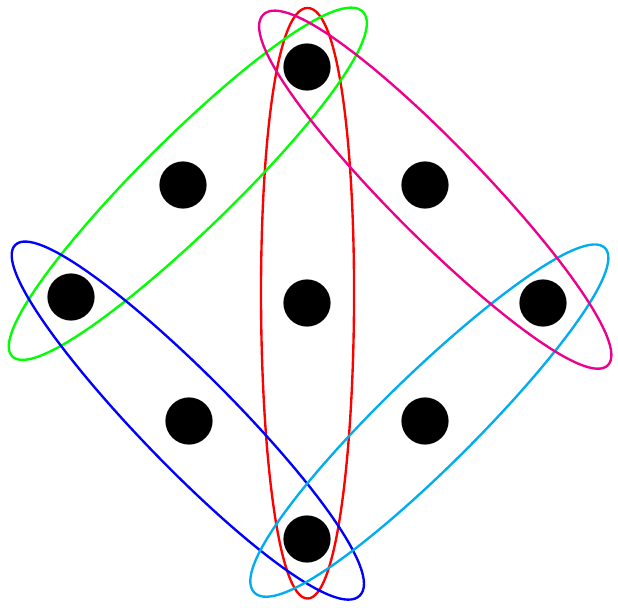}\par
\includegraphics[width=2cm]{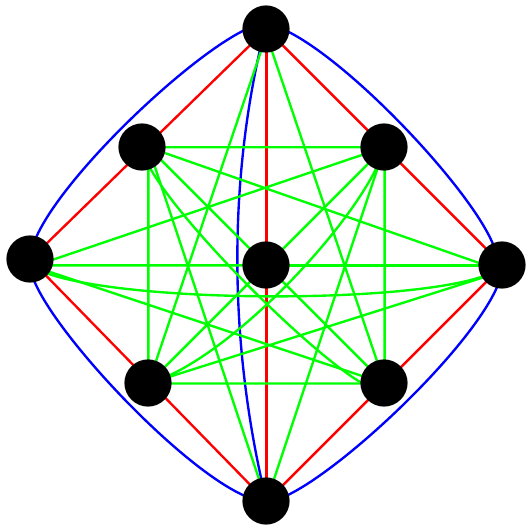}\par
\end{multicols}
\begin{multicols}{2}
    \begin{center}
        $H_7^2$ and $M_7^2$
    \end{center}\par
    \begin{center}
        $H_9$ and $M_9$
    \end{center}
\end{multicols}
\caption{These are the possible induced hypergraphs $H[X_i]$ of $H$ that respect the $f_0$-sparse condition, and their corresponding metric spaces. The metric spaces are labeled as $M_i$ or $M_i^j$, where \textcolor{red}{red} indicates distance 1, \textcolor{blue}{blue} distance 2 and \textcolor{green}{green} distance 1.5.}
\label{table1}
\end{figure}

\begin{theorem}[\textbf{Characterization of $\ceilhalf$-sparse 3-uniform hypergraphs}]
\label{thm:classification}
    Let $H = (V,E)$ be a $\ceilhalf$-sparse 3-uniform hypergraph. Then, for some $0\leq k\leq 2$ and for $i=1,\dots,k$, there are sets $X_i\subseteq V$ such that for any $i\neq j$, $|X_i\cap X_j|\leq 1$, and for each $i$, $H[X_i]$ is a copy of one of the hypergraphs $H_4$, $H_5^1$, $H_5^2$, $H_5^3$, $H_6$, $H_7^1$, $H_7^2$ or $H_9$ (presented in Figure \ref{table1}). Moreover, the hypergraph $H \setminus (\bigcup E[X_i])$ is $(6,2)$-sparse and any two of its hyperedges share at most one vertex. Also, if $k=2$ and $X_1\cap X_2\neq \emptyset$, then any edge of $H \setminus (\bigcup E[X_i])$ intersects $X_1\cup X_2$ in at most one vertex. 
\end{theorem}

Finally, we are able to prove our main result. The idea of the proof is, for a given $\ceilhalf$-sparse hypergraph $H=(V,E)$, we construct a metric $\rho$ using Theorem \ref{sixtwoNoGoat} and Figure \ref{table1} and prove that for the metric space $M=(V, \rho)$ we have $H_{M} = H$.

\begin{proof}[Proof of Theorem~\ref{thm:fsparse}]

If $H$ is $(6,2)$-sparse, then by Theorem \ref{thm:sixtwo} we know $H$ is metric. So, let us assume otherwise. By Theorem \ref{thm:classification}, there are some sets $X_1,\dots,X_k$ such that $H \setminus (\bigcup E[X_i])$ is $(6,2)$-sparse and two of its hyperedges share at most one vertex. Let $\rho_0$ denote the metric constructed for the hypergraph $H \setminus (\bigcup E[X_i])$ in the proof of Theorem \ref{sixtwoNoGoat}. The metric space $M_0 = (V, \rho_0)$ satisfies that $H_{M_0} = H \setminus (\bigcup E[X_i])$.

As seen in Figure \ref{table1}, $H[X_i]$ is metric for all $i=1,\dots,k$. Let $\rho_i$ denote the metric presented in Figure \ref{table1} in the corresponding metric space to $H[X_i]$. Now, we define the map $\rho: V \rightarrow \{ 0,1,1.5,2 \}$ as follows
\begin{equation*}
\rho(x,y) = 
    \begin{cases}
        \rho_i (x,y) &\text{if both $x,y$ are in $X_i$}, \\
        \rho_0 (x,y) &\text{otherwise.}
    \end{cases}
\end{equation*}
The map $\rho$ is well-defined since $H \setminus (\bigcup E[X_i])$ is $(6,2)$-sparse and $\rho_i$ are well-defined. It is also a metric. 
Any edge $\{ x,y,z \} \in E$ forms a collinear triple with respect to $\rho$ by Theorem \ref{thm:sixtwo} and Figure \ref{table1}. Now, assume $\{ x,y,z \} \notin E$. If $\{ x,y,z \} \subset X_i$ for some $i$ or if no two are in the same $X_i$, then it is clear that $\{ x,y,z \}$ is not a collinear triple. If this is not the case, Theorem~\ref{thm:classification} guarantees that at least one of the distances is $1.5$, so the triple is not collinear.
\end{proof}

\printbibliography

\end{document}